\documentclass[12pt,a4paper]{amsart}
\usepackage{amssymb,amsmath,amsthm}
\usepackage{graphicx,enumerate}
\usepackage{hyperref}
 \usepackage[color=green!40]{todonotes}

\newcommand\cC{{\mathcal C}}

\newcommand\cE{{\mathcal E}}

\newcommand\cS{{\mathcal S}}

\newcommand\cG{{\mathcal G}}
\newcommand\cH{{\mathcal H}}

\newcommand\cM{{\mathcal M}}

\newcommand{\ex}{\mathop{}\!\mathrm{ex}}

\makeatletter
\newtheorem*{rep@theorem}{\rep@title}
\newcommand{\newreptheorem}[2]{%
\newenvironment{rep#1}[1]{%
 \def\rep@title{#2 \ref{##1}}%
 \begin{rep@theorem}}%
 {\end{rep@theorem}}}
\makeatother

\newtheorem{thm}{Theorem}
\newtheorem{theorem}{Theorem}

\newtheorem{lemma}[thm]{Lemma}

\newtheorem{conjecture}[thm]{Conjecture}
\newtheorem{proposition}[thm]{Proposition}

\newtheorem{problem}[thm]{Problem}

\newtheorem{corollary}[thm]{Corollary}

\newtheorem{example}[thm]{Example}

\newtheorem{clm}[thm]{Claim}

\newcommand\cref[1]{Corollary~\ref{cor:#1}}

\title{Edge mappings of graphs: Tur\'an type parameters}

\author{Yair Caro}
\address{Department of Mathematics, University of Haifa-Oranim, Israel}
\email{yacaro@kvgeva.org.il}

\author{Bal\'azs Patk\'os}
\address{Alfr\'ed R\'enyi Institute of Mathematics}
\email{patkos@renyi.hu}
\thanks{Research of Patk\'os, Tuza, and Vizer is partially supported by NKFIH grants SNN 129364 and FK 132060.}

\author{Zsolt Tuza}
\address{Alfr\'ed R\'enyi Institute of Mathematics and University of Pannonia}
\email{tuza.zsolt@mik.uni-pannon.hu}

\author{M\'at\'e Vizer}
\address{Department of Computer Science and Information Theory, Budapest University of Technology and Economics}
\email{vizermate@gmail.com}

\date{}
\date{}

\begin{document}

\maketitle

\begin{abstract}
    In this paper, we address problems related to parameters concerning edge mappings of graphs. The quantity $h(n,G)$ is defined to be the maximum number of edges in an $n$-vertex graph $H$ such that there exists a mapping $f: E(H)\rightarrow E(H)$ with $f(e)\neq e$ for all $e\in E$ and further in all copies $G'$ of $G$ in $H$ there exists $e\in E(G')$ with $f(e)\in E(G')$. Among other results, we determine $h(n, G)$ when $G$ is a matching and $n$ is large enough.

    As a related concept, we say that $H$ is unavoidable for $G$ if for any mapping $f: E(H)\rightarrow E(H)$ with $f(e)\neq e$ there exists a copy $G'$ of $G$ in $H$ such that $f(e)\notin E(G')$ for all $e\in E(G)$. The set of minimal unavoidable graphs for $G$ is denoted by $\cM(G)$. We prove that if $F$ is a forest, then $\cM(F)$ is finite if and only if $F$ is a matching, and we conjecture that for all non-forest graphs $G$, the set $\cM(G)$ is infinite.

    Several other parameters are defined with basic results proved. Lots of open problems remain.
\end{abstract}

\section{Introduction}

Set-mappings or alternatively edge-mappings were first studied systematically by Erd\H os and Hajnal \cite{EH58} who were mostly interested in and motivated by mappings between infinite sets. Yet they also offered,  under the influence of Ramsey's theorem, the first instance of such problems in the finite case. Nevertheless,
while this line of research captured much attention and interest in the infinite case, only modest attention was given to the finite case. The latter was reconsidered more systematically in the 1980s \cite{AC, ACT89, C, CS},
but then the subject was again dormant for nearly 35 years (with few exceptions like \cite{AM06}) until recent interest in these types of problems was renewed by Conlon, Fox, and Sudakov \cite{CFS16,CFS20,CFS-x}.

Here, we start a more systematic study of problems about edge mappings in the context of graphs. The present manuscript focuses on parameters that are analogous to Tur\'an-type problems, and a parallel manuscript \cite{CPTV} addresses Ramsey-type parameters. 
We hope to raise various, simple to state yet perhaps not so simple to solve open problems regarding the parameters and functions that capture the essence of edge-mappings in parallel to those parameters and functions in Tur\'an problems.

\medskip

In order to present our main results we introduce the necessary definitions and notation.

\subsection{Definitions}
Let $E_n  = E(K_n)$  denote the edge set of the complete graph $K_n$ on $n$ vertices, and let $f: E_n \rightarrow  E_n$. A subgraph $G$ of $K_n$ is called 

\begin{itemize}

\item \textit{$f$-free,} if $f(e) \notin E(G)$ for all $e \in E(G)$, and
\item \textit{$f$-exclusive,} if $|f(e ) \cap   V(G) | = 0$ for all $e\in E(G)$.

\end{itemize}

 \medskip

\noindent
Let us introduce some notation for edge mappings. Let $F_n  :=   \{ f : E_n\rightarrow E_n \}$.  

\begin{itemize}

\item For $d = 0 ,1$ let $F_{n,d} := \{  f : E_n \rightarrow E_n  ~\text{such that} \ |f(e) \cap  e| \le d  ~\text{for every}\ e \in E_n \}$;

\smallskip 

\item for $d = 0 ,1$ and $0\le m \le \binom{n}{2}$ let 
$$F_{n,m,d}:=\{  f : E_n \rightarrow E_n  ~\text{such that for at least $m$ edges} \ e \in E_n , |f(e )\cap e|\le d   \};$$

\smallskip 

\item for $d = 0 ,1$ and a graph $H$ let $$F_{H,d}:=\{f:E(H)\rightarrow E(H) ~\text{such that for every} \ e \in E(H)  , \ |f(e )\cap e|\le d   \}.$$

\end{itemize}

\vskip 0.3truecm

Now we introduce further parameters related to Tur\'an numbers of graphs, that we will focus on during our investigations.  

\medskip 



First, let us recall that for $n \ge 1$ and a graph $F$, the \textit{Turán number} $\ex(n,F)$ is the maximum number of edges that an $n$-vertex graph can have without containing $F$ as a subgraph. 

\medskip 

We consider the following further extremal functions.

\begin{itemize}

\item For a given graph $G$, $d  = 0 ,1$ and  $n\ge 3$  let 
$$q(n,d,G) := \max\{ m : \exists \ f\in F_{n,m,d}, \ \text{with no $f$-free copy of G in $K_n$}  \}.$$  

We abbreviate $q(n,1,G)$ to $q(n,G)$.

\medskip  

\item For a given graph $G$, $d  = 0 ,1$ and $n \ge 3$  let
\begin{eqnarray}
h(n,d,G) := \max\{ m :& & \hspace*{-1.4em}
\text{$\exists$ $H$, $|V(H)|=n$, $|E(H)|=m$, and $\exists f \in F_{H,d}$} \nonumber \\
& &  \hspace*{-1.4em}
\text{ with no $f$-free copy of $G$ in $H$} \}. \nonumber
\end{eqnarray}

We abbreviate $h(n,1,G)$ to  $h(n,G)$.

\medskip  

\item 

For a given graph $G$, and $n \ge 4$   let 
\begin{eqnarray}
s(n, G) :=  \max\{ m :& &  \hspace*{-1.4em}
\text{$\exists$ $H$, $|V(H)|=n$, $|E(H)|=m$, and $\exists f\in F_{H,0}$ } \nonumber \\
& &  \hspace*{-1.4em}
\text{with no $f$-exclusive  copy of $G$ in $H$} \}. \nonumber
\end{eqnarray}

\end{itemize}

All these definitions naturally extend to families $\cG$ of graphs by requiring no $f$-free / \newline $f$-exclusive copies for any $G\in \cG$.

\medskip

\noindent \textbf{Further notation on graph operations}. For any pair $G, H$ of graphs, $G\cup H$ denotes their vertex-disjoint union, $G+ H$  their join, and for an integer $k$ we denote by $kG$ the vertex-disjoint union of $k$ copies of $G$.

\medskip 

In Sections \ref{h} and \ref{s}, we prove several upper and lower bounds on $h(n, G)$ and $s(n, G)$ for specific graphs. Here we state a theorem that determines the value of $h(n, G)$ for matchings.

\vbox{
\begin{thm}\label{hn2k2}
\
\begin{enumerate}

\item 
We have $h(n,2K_2)=n$ for $n \ge 7$, and for smaller vales of $n$ we have 
$h(4,2K_2)=6$ and
$h(5,2K_2)=h(6,K_2)=7$.

\item
For any $t\ge 3$ there exists $n_t$ such that if $n\ge n_t$, then $$h(n,tK_2)=\ex(n,tK_2)+t-1=(t-1)(n-t+1)+\binom{t-1}{2}+(t-1).$$
\end{enumerate}
\end{thm}
}

The following notion naturally arises in the study of $h(n,G)$. Given a graph $H$  we say that $G$  is \textit{$f$-free unavoidable for $H$} or \textit{unavoidable for $H$}, if for every $f \in F_{G,1}$, there exists an $f$-free copy of $H$. Observe that if $G$ is unavoidable for $H$  and $G$ is a subgraph of $Q$,  then $Q$ is also unavoidable for $H$, as for any $f \in F_{Q,1}$ one can restrict $f$ to a copy of $G$ and obtain an $f$-free copy of $H$. (If $f(e)\in E(Q)\setminus E(G)$ for some $e\in E(G)$, then it is even better, or formally we can then let $f(e)=e^*$ for some $e^*\neq e$, $e^*\in E(G)$.)
Therefore, it is enough to consider \textit{minimal} $H$-unavoidable graphs, where $G$ is minimal unavoidable for $H$ if it does not contain a proper subgraph $G^*$ which is unavoidable for $H$.

Let $\cM(H)$ denote the set of minimal unavoidable graphs for $H$.
 The next observation establishes a connection between the function $h$ and the Tur\'an number.  

\begin{proposition}\label{unavoid}
    For any graph $H$ we have $h(n,H) = \ex(n, \cM(H) )$.
\end{proposition} 

We determine the families $\cM(P_3)$ and $\cM(2K_2)$.
In the general setting, our main result on $\cM(G)$ is the following theorem.

\begin{thm}\label{forest}
    If $F$ is a forest, then $\cM(F)$ is finite if and only if $F = tK_2$.
\end{thm}

We prove Theorem \ref{forest} in Section \ref{m(G)} along with other statements on $\cM(G)$ for non-forest $G$. We believe that Theorem \ref{forest} is valid not only for forests.

\begin{conjecture}\label{mconj}
$\cM(G)$ is finite if and only if $G = tK_2$.
\end{conjecture}

The remainder of the paper is structured as follows. In Section \ref{prelim}, we establish some general connections between the parameters of interest, in particular, we prove Proposition \ref{unavoid}. Proofs of results on $h(n, G)$, $s(n, G)$ and $\cM(G)$ are contained in Sections \ref{h}, \ref{s}, and \ref{m(G)}, respectively, while Section \ref{comm} is devoted to some concluding remarks and open problems.

\section{Preliminaries and general connections between parameters}\label{prelim}

Alon and Caro \cite{AC} proved the following sandwich theorem for Tur\'an numbers and the function $h$.

\begin{theorem}[Alon, Caro \cite{AC}]\label{acex}
For any family $\cG$ of graphs, we have $$\ex(n,\cG)\le h(n,\cG)\le 3\cdot \ex(n,\cG).$$ Moreover, for graphs $G$ with $\chi(G)\ge 3$, we have $h(n,G)=(1+o(1))\ex(n,G)$.
\end{theorem}

Next, we examine whether equality is possible in the inequalities $\ex(n,\cG)\le h(n,\cG)\le 3\cdot \ex(n,\cG)$. The following example shows that the upper bound in Theorem \ref{acex} is attainable. 

\begin{example}
Motivated by the term of ``delta-systems'' we use the notation $\Delta(r,t)  = \{  K_{1,r}, tK_2\}$.  It is known \cite{CH} that $\ex(n,\Delta(r,t))=(r-1)(t-1)$ if $r\ge 2t+1$. We claim that
$h(n,\Delta(r,t))=3(r-1)(t-1)$ holds. As the upper bound is given by 
 Theorem \ref{acex}, we only need to see the lower bound. We take  $G=(t-1)K_{1,3(r-1)}$. As it does not contain $tK_2$, we only need to define a mapping $f\in F_{G,1}$ with no $f$-free $K_{1,r}$. If in the $i$th star the edges are $e_{3i-2},e_{3i-1},e_{3i}$, then for any $i=1,2,\dots,r-1$ and $j=1,2$ we let $f(e_{3i-j})=e_{3i-j+1}$ and $f(e_{3i})=e_{3i-2}$.
\end{example}

On the other hand, the next proposition shows that the lower bound in Theorem \ref{acex} is never attained (except for the trivial case of $G=K_2$).

\begin{proposition}\label{strict}
For any graph $G$ with at least two edges, we have $h(n,G)>\ex(n,G)$ for all $n\ge |V(G)|$.
\end{proposition}

\begin{proof}
Let $H$ be an $n$-vertex $G$-free graph with $ex(n,G)$ edges, and let $H_{e^*}=H\cup \{e^*\}$ for any $e^*\notin E(H)$. Finally, set $f(e)=e^*$ for every $e\in E(H)$ and $f(e^*)=e'$ for some $e'\in E(H)$. Clearly $f\in F_{H_{e^*},1}$, and by the $G$-free property of $H$, any copy of $G$ in $H_{e^*}$ must contain $e^*$, but as there exists another edge in $G$, it cannot be $f$-free.
\end{proof}

\begin{proposition}\label{qh}
For any graph $G$ we have $q(n, G) = h(n, G)$.
\end{proposition}

\begin{proof}
 Clearly by definition $q(n, G) \ge h(n, G)$ as a function realizing $h(n, G)$ is a function allowed for  $q(n, G)$. 

Suppose $f$ is a function realizing $q(n, G)$ and let $H$ be the graph whose edges are those for which $f( e) \neq e$.  Suppose for some $e \in  E(H)$ we have $f(e) = e_1\notin E(H)$.  Choose an edge $e^* \in E(H) \setminus \{ e\}$  and change $f$ to $f^*$ by defining  $f^*(e):=f(e)$ for all edges of $H$ except  that $f(e)=e^*$. 

If there is an $f^*$-free copy of $G$ in $H$, then this is an $f$-free copy of $G$ as well, as $e_1$ is not in $E(H)$ hence not in $E(G)$. So we can replace all edges that are mapped out of $E(H)$  to non-fixed edges that are mapped into $E(H)$ not creating $f$-free copies of $G$.  
\end{proof}

\begin{proof}[\textbf{Proof of Proposition \ref{unavoid}}]
    Let $G$ be any graph realizing $h(n, H)$.  Then $G$ cannot contain any member of $\cM(H)$, because if $G$ contains $G^*\in \cM(H)$, then  $G$ is also unavoidable for $H$ and so for any $f\in F_{G,1}$ there exists an $f$-free copy of $H$. Hence $h(n,H)  \le \ex(n,\cM(H))$.

On the other hand, let $G$ be any $\cM(H)$-free graph having $\ex(n,\cM(H))$ edges. Then  $G$ contains no member of $\cM(H)$, and so it is not unavoidable for $H$, hence there exists an edge mapping  $f\in F_{G,1}$ with no $f$-free copy of $H$.
Thus $h(n,H)  \ge \ex(n,\cM(H))$.
\end{proof} 

We finish this section with a simple lemma. Part (1) is from \cite{AC} (Lemma 2.3), but as we will use it multiple times in later proofs, we restate and reprove it in order to be self-contained.
The undirected underlying graph of a directed graph $\overrightarrow{G}=(V,\overrightarrow{E})$ is $G=(V,E)$, where $E=\{(uv): \overrightarrow{uv} \in \overrightarrow{E} ~\text{or}\ \overrightarrow{vu}\in\overrightarrow{E}\}$.


\begin{lemma}\label{outdeg}
    Let $\overrightarrow{\Gamma}$ be a directed graph on $n$ vertices with maximum out-degree at most $d$ and let $\Gamma$ be the undirected underlying graph of $\overrightarrow{\Gamma}$.
\begin{enumerate}
    \item
    Then we have $\chi(\Gamma)\le 2d+1$, and thus $\alpha(\Gamma)\ge \frac{n}{2d+1}$.
    \item
    If $d=1$, and $m$ vertices have out-degree $0$, then  $\alpha(\Gamma)\ge m+\lceil \frac{n-2m}{3}\rceil$. Furthermore, if $m=0$ and $3|n$, then equality holds if and only if $\overrightarrow{\Gamma}$ consists of vertex-disjoint cyclically directed triangles.
    \end{enumerate}
\end{lemma}

\begin{proof}
    We start with the proof of (2). Note that as the out-degree of each vertex of  $\overrightarrow{\Gamma}$ is at most 1, the components of $\Gamma$ are trees or unicyclic\footnote{A unicyclic graph is a connected graph containing exactly one cycle.} graphs. Also, note that the components containing a vertex with 0 out-degree (in $\overrightarrow{\Gamma}$) are trees. Trees are bipartite, unicyclic graphs are 3-colorable. This implies inequality. 
    
    If $m=0$, then we need all components of $\Gamma$ to be unicyclic. As all such graphs can be three-colored with one color class having size one, we need $\left\lceil \frac{|C|-1}{2} \right\rceil=\frac{|C|}{3}$ for all components $C$. This implies $|C|=3$ and the out-degree condition yields that the triangles are cyclically directed.

\medskip    

    To see (1), observe that if the out-degree of $\overrightarrow{\Gamma}$ is $d$, then the average degree, and thus the degeneracy\footnote{ A $k$-degenerate graph is an undirected graph in which every subgraph has a vertex of degree at most $k$. The degeneracy of a graph is the smallest value of $k$ for which it is $k$-degenerate.} of $\Gamma$ is at most $2d$, so its chromatic number is at most $2d+1$.
\end{proof}


We will apply Lemma \ref{outdeg} in the following situation: given are a graph $G$, a mapping $f\in F_{G,1}$, and a partition $\cE=\{E_1,E_2,\dots,E_m\}$ of $E(G)$. We define $\overrightarrow{\Gamma}_{f,\cE}$ to have vertex set $\cE$ and there is a directed edge from $E_i$ to $E_j$ if and only if there exists $e\in E_i$ with $f(e_i)\in E_j$. If all $E_i$s consist of a single edge, then we will drop $\cE$ from the subscript and write $\overrightarrow{\Gamma}_f$. An edge subset $E'\subset E(G)$ is $f$-free if and only if the corresponding vertices form an independent set in $\overrightarrow{\Gamma}_f$.


\section{Results about $h(n,G)$}\label{h}

In this section, we obtain results on $h(n, G)$. Let us mention that prior to our paper the only graph class for which the exact $h$-value was known is that of stars due to Alon and Caro \cite{AC}. They showed $h(n,K_{1,r}) = \min \{ \binom{n}{2}  , n(r-1) \}$ for all $r\ge 1$. We start our investigations with general observations strengthening the lower bound of Proposition \ref{strict}.

\begin{proposition}
\label{p:lba}

\ 

\begin{enumerate}

\item \label{indmatch} Let $G$ be a graph with no isolated edges. Then $$h(n,G)\ge \ex(n,G)+\lfloor\alpha(H)/2\rfloor$$ for any $G$-free graph $H$ on $n$ vertices with $\ex(n,G)$ edges.

\item \label{starexp}
Let $T$ be a non-star tree and $H$ be a $T$-free graph on $k$ vertices. Then $$h(n,T)\ge \left\lfloor\frac{n-1}{k} \right\rfloor(e(H)+\alpha(H)).$$

\item \label{cor}
For any tree $T$ on $k$ vertices, we have $$h(n,T)\ge \Big(\frac{k-2}{2}+ \frac{1}{k-1}\Big)n-c_k$$ for some constant $c_k$ that depends only on $k$ and not on $T$.

\end{enumerate}

\end{proposition}

\begin{proof}
To prove (1) let $H$ be an $n$-vertex $G$-free graph with $ex(n,G)$ edges, and let $I$ be an independent set in $H$ of size $\alpha(H)$. Then pick any matching $M$ that covers $2\lfloor \alpha(H)/2\rfloor$ vertices of $I$. Clearly, $M\cap E(H)=\emptyset$. Define $H'$ by $V(H')=V(H)$, $E(H')=E(H)\cup M$ and set $f\in F_{n,H',1}$ as follows: if $e\in E(H)$ is incident to some edge $e^*\in M$, then let $f(e^*)=e$, for other edges $e$ of $E(H)$ or for edges of $M$ define $f$ arbitrarily but satisfying $f(e)\neq e$. Then any copy of $G$ in $H'$ must contain at least one edge $e^* \in M$ and, as $e^*$ is not an isolated edge of $G$, an edge $e$ incident to $e^*$. As $f(e)=e^*$, the copy of $G$ cannot be $f$-free. 

\medskip 

To prove (2) let $I$ be an independent set of $H$ of size $\alpha(H)$. Define $H'$ to be a graph that contains a special vertex $u$, and $H'\setminus \{u\}$ having $\lfloor \frac{n-1}{k}\rfloor$ copies of $H$ possibly with some isolated vertices. The special vertex $u$ is joined to the vertices in all copies of $I$. For any edge $e=xy$ with $x$ in $I$ in some component and $y\neq u$, we set $f(xy)=xu$ and for all other edges $e$ of $H'$, we define $f$ such that $f(e)\neq e$ holds. As $H$ is $T$-free, any copy of $T$ must contain an edge $xu$. Also, as $T$ is not a star, it must contain an edge that is incident to an edge not containing $u$. So we may assume that the copy of $T$ must contain $xy$ for some $y\neq u$, and thus the copy of $T$ cannot be $f$-free.

\medskip 

Part (3) follows from (\ref{starexp}) with $H=K_{k-1}$.
\end{proof}





For a tree $T$ with bipartition $(A,B)$, we set $\Delta(A)=\max\{d_T(u):u\in A\}$, $\Delta(B)=\max\{d_T(u):u\in B\}$, and $\Delta^*(T)=\min\{\Delta(A),\Delta(B)\}$. With these definitions, we obtain the following constructive lower bounds.

\begin{proposition}\label{delta*}

\ 

\begin{enumerate}

\item \label{maxdeg}
For any tree $T$ and $n> 2(\Delta(T)-1)$, we have $h(n,T)\ge n(\Delta(T)-1)$.

\item \label{minmax}
For any tree $T$, we have $h(n,T)\ge (3\Delta^*(T)-3)(n-3\Delta^*(T)+3)$. 

\end{enumerate}

\end{proposition}

\begin{proof}
To prove (1) it is known that there exists a $2(\Delta(T)-1)$-regular graph $H$ on $n$ vertices. Consider a closed Eulerian trail $e_1,e_2,\dots,e_{n(\Delta(T)-1)}$ of $H$ and define $f(e_i)=e_{i+1}$ implying $f\in F_{n,H,1}$. If $u$ is a maximum-degree vertex of $T$, then consider the vertex corresponding to $u$ in any copy of $T$ in $H$. As there are exactly $\Delta(T)-1$ incoming and outgoing edges in the trail at $u$, there must exist a pair $e_j,e_{j+1}$ that belongs to the copy of $T$. Therefore this copy cannot be $f$-free.

\medskip 

To prove (2) consider $K=K_{3(\Delta^*(T)-1),n-3(\Delta^*(T)-1)}$ and name the vertices in the part of size $3(\Delta^*(T)-1)$ as $a_i,b_i,c_i$ for $i=1,2,\dots,3(\Delta^*(T)-1)$. Define the mapping $f$ by $f(a_ix)=b_ix$, $f(b_ix)=c_ix$, and $f(c_ix)=a_ix$ for all $i=1,2\dots,3(\Delta^*(T)-1)$ and all vertices $x$ being in the other part of $K$. Any copy of $T$ must contain a vertex $x$ in the big part of $K$ of degree $\Delta^*(T)$. Therefore, there must exist an index $i$ such that at least two of $a_ix,b_ix,c_ix$ belong to the copy of $T$, so this copy cannot be $f$-free.
\end{proof}




\begin{corollary}\label{doublestar}
For the balanced double star $D_{k,k}$, we have $3kn-9k^2\le h(n, D_{k,k})\le 3kn$.
\end{corollary}

\begin{proof}
To obtain the lower bound, apply Proposition \ref{delta*} (\ref{minmax}) with $\Delta^*(D_{k,k})=k+1$. The upper bound follows from Theorem \ref{acex}, as for double-stars the Erd\H os--S\'os conjecture is known to hold.
\end{proof}

\begin{proposition}\label{paths}
We have
\begin{enumerate}
    \item 
    $h(n,P_4)\ge 3(n-3)$,
    \item
    $h(n,P_7)\ge h(n,P_6)\ge h(n,P_5)\ge 3n-6$,
    \item
    $h(n,P_{2m+1})\ge h(n,P_{2m}) \ge m(n-m)$ for all $m\ge 4$.
\end{enumerate}
\end{proposition}

\begin{proof}
$P_4=D_{1,1}$, so Corollary \ref{doublestar} yields (1).

\medskip 

Also, (3) follows from Proposition \ref{p:lba} (\ref{starexp}) with $H=K_{m-1,n-m}$, $k=n-1$, and $\alpha(K_{m-1,n-m})=n-m$. 
(Actually, the graph verifying the lower bound is $H=K_{m,n-m}$.)

\medskip 

Finally to see (2), consider $K_3+E_{n-3}$ with the vertices of the $K_3$ being $a,b,c$. Define the mapping $f$ by $f(ab)=bc$, $f(bc)=ca$, $f(ca)=ab$ and for any $x$ in $V(E_{n-3})$ we let $f(ax)=bx$, $f(bx)=cx$, $f(cx)=ax$. Then if a copy of $P_5$ uses two edges of $K_3$, then it is not $f$-free. If such a copy uses at most one such edge, then there exists an $x\in V(E_{n-3})$ such that $x$ is incident to two edges of the $P_5$, and thus this copy of $P_5$ cannot be $f$-free either. 
\end{proof}

Observe that if $n$ is not divisible by 3, then $\ex(n,P_4)=n-1$ and so by Theorem \ref{acex}, we have $h(n,P_4)\le 3n-3$.


\medskip

In the remainder of this section, we prove Theorem \ref{hn2k2}. We repeat the statements of both of its parts, separately.

\begin{thm}

We have $$h(n,2K_2)=n$$ for $n \ge 7$.
Furthermore, $h(4,2K_2)=6$, and $h(5,2K_2)=h(6,K_2)=7$.

\end{thm}

\begin{proof}
The lower bound $h(4,2K_2)\ge 6$ is shown by the mapping $f:E(K_4)\rightarrow E(K_4)$ with $f(e)=V(K_4)\setminus e$. To see $h(6,2K_2)\ge h(5,2K_2)\ge 7$, consider the graph $H$ that we obtain from $K_4$ by adding a pendant edge $e$ and the mapping $f\in F_{H,1}$ with $f(e^*)=e$ for all $e^*\in E(H)$ disjoint with $e$, and $f(e^{**})=V(K_4)\setminus e^{**}$ for all other $e^{**}\in E(K_4)$ (and $f(e)$ defined arbitrarily). Finally, the general lower bound $h(n,2K_2)\ge n$ is given by the graph $H$ obtained from $K_{1,n-1}$ by adding an edge between two leaves and the mapping that maps all the edges of $K_{1,n-1}$ to the additional edge.

The upper bound $h(4,2K_2)\le 6$ is trivial as there are at most 6 edges on 4 vertices.
For larger $n$ we apply the following assertion whose validity is seen by observing that each pair $(e,f(e))$ can make only one copy of $2K_2$ non-$f$-free.

\begin{clm}\label{2k2claim}
    For a graph $H$ let $m$ denote the number of copies of $2K_2$ and $m'$ the number of edges that are contained in at least one copy of $2K_2$. If $m>m'$, then $H$ is unavoidable for $2K_2$.
\end{clm}

Based on this claim we first derive $h(5,2K_2)\leq 7$.
Observe that $K_5$ contains 15 copies of $2K_2$ because each of its edges occurs in precisely three copies.
Hence a graph with 5 vertices and 8 edges contains at least 9 copies, more than the number of edges.

For $n\geq 6$ we prove the upper bounds by induction.
Let us prove first $h(6,2K_2)\leq 7$.
Should it be 8, then a graph verifying it would have no isolated vertices, due to $h(5,2K_2)=7$.

The 8 edges determine 28 edge pairs, and a degree sequence $d_1\geq \cdots \geq d_6$ yields $s:=\sum \binom{d_i}{2}$ pairs that are not $2K_2$.
To have a chance without $f$-free copies of $2K_2$ we would need $s\geq 20$ because of Claim \ref{2k2claim}.
However, the largest degree is at most 5 and the smallest one is at least 1.
The sequence cannot begin with 5, 5 because it would already require 9 edges.
Therefore the degree sequence with maximum $s$ is $5,4,2,2,2,1$ whose $s$-value is 19 only (the sequence $5,4,3,2,1,1$ is not graphical).

We prove $h(7,2K_2)\leq 7$ in a similar way.
The degree sequence cannot begin with $6,4$ as it would require 9 edges. The sequence $6,3,3,1,1,1,1$ is not graphical.
So the sequence is $6,3,2,2,1,1,1$ with $s=20$, or $6,2,2,2,2,1,1$ ($s=19$) or $5,4,2,2,1,1,1$ ($s=18$) or something with an even smaller $s$.
Only the first sequence has a chance, but the edge $e$ joining the first two vertices is not contained in a $2K_2$, hence the pair $(e,f(e))$ is useless and an $f$-free $2K_2$ occurs.

The induction step for $n>7$ is really simple because the formula would become false only if a graph with minimum degree of 2 would be found.
However, with $n+1$ edges, the degree sum is only $2n+2$, hence $s$ would be maximized with the degree sequence $4,2,2,\dots,2$.
The corresponding $s$ is only $n+5$ while we have just $n+1$ edges, however $n+1 < \binom{n+1}{2} - (n+5)$ holds for every $n \geq 6$.
This completes the proof.
\end{proof}

Finally, we restate and prove the second part of Theorem \ref{hn2k2}.

\begin{thm}
For any $t\ge 3$ there exists $n_t$ such that if $n\ge n_t$, then $$h(n,tK_2)=\ex(n,tK_2)+t-1=(t-1)(n-t+1)+\binom{t-1}{2}+(t-1).$$
\end{thm}

\begin{proof}
We start with the general lower bound: consider the split graph $K_{t-1} + E_{n- t+1}$. Let the vertices of $K_{t-1}$ be $u_1,u_2,\dots, u_{t-1}$.  Embed a star on $t-1$ edges in $E_{n-t +1}$  with center $v$ and $t-1$ leaves $v_1,v_2,\dots,v_{t -1}$. Call this graph $H$. We define the mapping $f$ as follows: for any $1\le i\le t-1$ and $x\neq u_1,u_2,\dots,u_{t-1}$, we set $f(u_ix)=vv_i$ and define $f$ on all other edges $e$ with $f(e)\neq e$ arbitrarily. Any $tK_2$ in $H$ contains exactly one edge incident with $u_1,u_2,\dots,u_{t-1},v$ each and these edges must be distinct. So if the edge incident with $v$ is $vv_i$, then for the edge $e_i$ incident with $u_i$ we have $f(e_i)=vv_i$ and thus this copy is not $f$-free. 

To see the general upper bound, we first prove a claim.

\begin{clm}\label{killingmatching}
For any $t\ge 2$ there exists a constant $\alpha_t$ such that for any $n$-vertex graph $H$ with $\ex(n,tK_2)$ edges, $f\in F_{H,1}$ and $e\in E(H)$ there are at most $\alpha_tn^{t-2}$ copies of $tK_2$ that contain both $e$ and $f(e)$.
\end{clm}

\begin{proof}[Proof of Claim]
The $t$-matching should contain $t-2$ other edges of $H$, and $\ex(n,tK_2)$ is a linear function of $n$ for every fixed $t$.
\end{proof}

Now we are ready to prove the general upper bound of the theorem. Let $H$ be a graph on $n$ vertices with $\ex(n,tK_2)+t$ edges. Then for any $T\subset V(H)$ with $|T|=t-1$, there exist at least $t$ edges not incident with any vertex of $T$.

Let $d_1\ge d_2\ge \dots \ge d_n$ be the degree sequence of $H$ and let $u_1,u_2,\dots,u_n$ be the corresponding vertices. Fix $f\in F_{H,1}$ and consider the following cases.

\vskip 0.3truecm

Case 1. $\frac{2}{7}n\le d_{t-1}\le n-D_t$ for some appropriately chosen constant $D_t$.

\vskip 0.15truecm

Any edge $e$ that is not incident with any $u_i$ is contained in $\prod_{j=1}^{t-1}(d_j-2t)$ copies of $tK_2$ such that $e$ is the only such edge in the $t$-matching. We have at least $D_t$ such edges, so $H$ contains at least $D_tn^{t-1}/4^{t-1}$ copies of $tK_2$. By Claim \ref{killingmatching}, every edge can ruin the $f$-free property of at most $\alpha_tn^{t-2}$ copies of $tK_2$, so the number of non-$f$-free copies of $tK_2$ is at most $t\alpha_tn^{t-1}$ which is smaller than the total number of copies if $D_t\ge t\cdot 4^{t-1}\cdot\alpha_t$. 

\vskip 0.3truecm

Case 2. $d_{t-1}\ge n-D_t$.

\vskip 0.15truecm

Let us fix edges $e_1,e_2,\dots,e_t$ none of which is incident with any of $u_1,u_2,\dots,u_{t-1}$. Let $E^0_i=\{u_ix:x\neq u_j ~j\neq i, ~x\notin e_h\}$, so $|E^0_i|\ge n-D_t-3t$. Observe that for any $i$ there is at most one $j=j_i$ such that for at least $n/2$ edges $e\in E^0_i$ we have $f(e)=e_j$. Therefore there exists a $j^*$ that is not a $j_i$ for any $1\le i \le t-1$. Then we set $e_{j^*}$ to be the first edge of our future $f$-free $tK_2$ and obtain $E^1_i$ from $E^0_i$ by removing $f(e_{j^*})$ and all edges $e$ with $f(e)=e_{j^*}$. By definition, $|E^1_i|\ge n/2-D_t-4t$. Next, for any $1\le i\le t-2$ there exists at most one edge $e^1_i\in E^1_{t-1}$ such that $f(e)=e^1_i$ for more than $n/4$ edges $e\in E^1_i$. Pick an edge $e^{t-1}\in E^1_{t-1}$ that is not $e^1_i$ for any $1\le i\le t-2$, and $e^{t-2}$ is our next edge of the future $tK_2$. We obtain $E^2_i$ from $E^1_i$ by removing $f(e^{t-2})$ and all edges $e$ for which we have $f(e)=e^{t-2}$. By definition, we have $|E^2_i|\ge n/4-D_t-5t$. We keep repeating this process to obtain an $f$-free copy of $tK_2$.

\vskip 0.3truecm

Case 3. $d_{t-1}\le \frac{n-t}{3}$

\vskip 0.15truecm

We will need the following claim.

\begin{clm}\label{supersatmatching}
For any $t\ge 2$ there exists a constant $C_t$ such that any $n$-vertex graph $H$ with $\ex(n,tK_2)+1$ edges and maximum degree at most $n/3$ contains at least $C_t\cdot n^{t}$ copies of $tK_2$.
\end{clm}

\begin{proof}[Proof of Claim]
We proceed by induction on $t$. If $t=2$  then any edge can be extended to a $2K_2$ with at least $|E(H)|-2\cdot n/3$ other edges, so the number of copies of $2K_2$ is quadratic.

For the inductive step, we argue similarly. For any edge $uv$ the graph $H\setminus \{u,v\}$ has more than $\ex(n-2,(t-1)K_2)$ edges, and thus by induction more than $C_{t-1}(n-2)^{t-1}$ copies of $(t-1)K_2$ each of which can be extended to a $tK_2$ with the edge $uv$. So any edge is extendable to at least $C_{t-1}(n-2)^{t-1}$ copies of $tK_2$, therefore $H$ contains at least $C_{t-1}(n-2)^{t-1}|E(H)|/t\ge C_tn^{t}$ copies of $tK_2$.
\end{proof}

Let $j$ be the smallest index such that $d_j\le \frac{n-t}{3}$, so $j\le t-1$. Then $H':=H-\{u_1,u_2,\dots,u_{j-1}\}$ has more than $\ex(n-(j-1),(t-(j-1))K_2)$ edges. So by Claim \ref{supersatmatching}, $H'$ contains at least $C_{t-(j-1)}\cdot (n-j+1)^{t-(j-1)}$ copies of $(t-(j-1))K_2$. Each such copy can be extended to at least $(n/3-3t)^{j-1}$ copies of $tK_2$ greedily with edges incident with $u_{j-1},u_{j-2},\dots,u_1$. Therefore $H$ contains at least $C'_tn^t$ copies of $tK_2$, so by Claim \ref{killingmatching}, there must exist an $f$-free copy of $tK_2$ if $n$ is large enough. 
\end{proof}

\section{Results about $s(n,G)$.}\label{s}

 Our first result asymptotically determines $s(n, G)$  for non-bipartite graphs, along the lines of how the asymptotic behavior of $h(n, G)$ is determined in \cite{AC}.

\begin{thm}
    For any graph $G$ with $\chi(G)\ge 3$, we have $s(n,G)=(1+o(1))\ex(n,G)$. If $G$ is bipartite, then $s(n,G)=o(n^2)$.
\end{thm}

\begin{proof}
    The inequality $\ex(n, G)\le s(n, G)$ holds by definition. To see the upper bound, assume $H$ is an $n$-vertex graph with $(1+\varepsilon)\ex(n,G)$ edges, and $f\in F_{H,0}$. By the celebrated result of Erd\H os and Stone \cite{ESt}, for any positive integer $q$, there exists $n_0=n_0(q)$ such that if $n\ge n_0$, then $H$ contains a copy of the complete $\chi(G)$-partite graph $K=K_{q,q,\dots,q}$ with $q$ vertices in each part. Observe that $K$ contains $\Theta(q^{|V(G)|})$ copies of $G$. On the other hand, for every edge $e\in E(K)$ there exist $O(q^{|V(G)|-3})$ copies of $G$ in $K$ that contain $e$ and that meet $f(e)$. Therefore there are $O(q^{|V(G)|-1})$ copies of $G$ in $K$ that are \textit{not} $f$-exclusive, so there must exist at least one $f$-exclusive copy of $G$.

    The statement about bipartite $G$ follows along the same lines as $ex(n,K_{q,q})=o(n^2)$ for any fixed $q$.
\end{proof}

\begin{proposition}\label{exclstar}
    Let $G$ be a graph in which no edge is incident to all other edges. If $\Delta(G)\ge 5r-4$, then any $f\in F_{G,0}$ admits an $f$-exclusive copy of $K_{1,r}$.
\end{proposition}

\begin{proof}
    Let $v$ be a vertex of $G$ with $d(v)\ge 5r-4$ and let $e_1,e_2,\dots,e_t$ be the edges incident to $v$. We will apply Lemma \ref{outdeg} to the directed graph $\overrightarrow{\Gamma}$ with $V(\overrightarrow{\Gamma})=\{e_i:i=1,2,\dots,t\}$ with $e_ie_j$ being an arc if and only if $f(e_i)\cap e_j\neq \emptyset$. As $f\in F_{G,0}$, $v\notin f(e_i)$ for all $i$, so the maximum out-degree of $\overrightarrow{\Gamma}$ is at most 2. Lemma \ref{outdeg} implies the existence of an independent set in $\Gamma$ of size $\lceil\frac{t}{5}\rceil\ge r$ that corresponds to an $f$-exclusive star of size $r$ with center $v$.
\end{proof}

If a graph $G$ on $n$ vertices contains an edge adjacent to all other edges, then $e(G)\le 2n-3$, so the following is an immediate corollary of Proposition \ref{exclstar}.

\begin{corollary}\label{corstar}
\

   $s(n,K_{1,r})\le \min \{ \binom{n}{2}  , 5(r-1)n/2 \}$.
   
\end{corollary}

\begin{proposition}
    \
    $s(n,K_{1,2})=\frac{5n}{2}$ for every $n$ divisible by 6. Moreover, for any $j=1,\dots,5$ there exists a constant $c_j$ such that $s(6n+j,K_{1,2})=15n+c_j$ if $n$ is large enough.
\end{proposition}

\begin{proof}
    The upper bound follows from Corollary \ref{corstar}. The lower bound is given by $n/6$ copies of $K_6$ with an $f$ defined as follows: $E(K_6)$ decomposes into five copies of $3K_2$. Within each $3K_2$, we map the edges cyclically. Each $3K_2$ destroys 12 distinct copies of $K_{1,2}$,  so the five copies destroy all  60 copies of $K_{1,2}$ in $K_6$. 
    
    The general lower bound follows from $s(6n+j,K_{1,2})\ge s(6(n-1),K_{1,2})+s(6,K_{1,2})=s(6(n-1),K_{1,2})+15.$ As Corollary \ref{corstar} implies $s(6n+j,K_{1,2})\le 15n+15$ there exist only a finite number of values of $n$ for which $s(6n+j,K_{1,2})-15n>s(6(n-1)+j,K_{1,2})-15(n-1)$.
\end{proof}

\begin{proposition}\label{texcl} For $r,t \ge 1$ we have
    $s(n, tK_{1,r} )  \le  \ex(n, ((4r +1)(t-1)+1)K_{1,5r-4})$. In particular, $s(n,tK_2)  \le \ex(n,(5t-4)K_2) .$
\end{proposition}

\begin{proof}
    Suppose $G$ is a graph on $n$ vertices with at least $\ex(n, ((4r +1)(t-1)+1)K_{1,5r-4}) + 1$ edges, and let $f\in F_{G,0}$. Then $G$ contains $(4r +1)(t-1)+1$ vertex-disjoint copies of $K_{1,5r -4}$.  The proof of Proposition \ref{exclstar} shows that each copy of $K_{1,5r-4}$ contains an $f$-exclusive copy of $K_{1,r}$.

Define a directed graph $\overrightarrow{\Gamma}$ where a vertex corresponds to a copy of $K_{1,r}$ and there is an $ij$ arc if there is an edge $e$  in the $i$-th copy of $K_{1,r}$  such that $f( e )$  contains at least one vertex from the $j$th copy of $K_{1,r}$. The outdegree of every vertex in $\overrightarrow{\Gamma}$ is at most $2r$. 
Hence, by Lemma \ref{outdeg}, we obtain a set of at least $t$ independent vertices corresponding to an $f$-exclusive copy of $tK_{1,r}$.
\end{proof}

Before our last result concerning the parameter $s$, let us recall that we use the notation $\Delta(r,t)  = \{  K_{1,r}, tK_2\}$.

\begin{proposition}
    $s(n,\Delta(r,r))=\Theta(r^2)$.
\end{proposition}

\begin{proof}
    Abbott, Hanson, and Sauer \cite{AHS} proved that $\ex(n,\Delta(r,r))=\Theta(r^2)$ (they determined the exact value of the extremal function). So if $G$ has more edges than $\ex(n,\Delta(5r-4,5r-4))$, then $G$ contains either a $K_{1,5r-4}$ or a $(5r-4)K_2$. In the former case Proposition \ref{exclstar}, in the latter case the proof of Proposition \ref{texcl} yields an $f$-exclusive $K_{1,r}$ or $rK_2$ for any $f\in F_{G,0}$.

    Clearly, $\ex(n,\Delta(r,r))$ is a lower bound on $s(n,\Delta(r,r))$.
\end{proof}

\section{Results about $\cM(G)$}\label{m(G)}

The notion of $\cM(G)$ is introduced and discussed in this research for the first time. 
Our main motivation in this section is to prove Theorem \ref{forest}. However, we start with two results that completely determine $\cM(G)$ for the two simplest graphs $P_3$ and $2K_2$.


\begin{proposition}
    $\cM(P_3)= \{  K_{1,4}\}\cup\{  C_k^+:k\ge 3\}\cup \{ B_k: k \ge 1  \}$, where $C_k^+$ is the cycle of length $k$ with a pendant edge and $B_k$ is the subdivision of the double-star $D_{2,2}$ with the degree-3 vertices being distance $k$ apart.
\end{proposition}

\begin{proof}
Clearly, any $G\in \cM(P_3)$ is connected. Also, if a graph $G$ contains a vertex $v$ with degree of at least $4$, then $v$ will be the middle vertex of an $f$-free copy of $P_3$ for any $f\in F_{G,1}$. As one can define an edge mapping of $K_3$ without $f$-free copies of $P_3$ by $e_1 \to e_2 \to e_3 \to e_1$, we have that $K_{1,4}$ is the only member of $\cM(P_3)$ containing a vertex of degree 4.

Next, we show that subgraphs of $C_k^+$ and $B_k$ are avoidable for $P_3$. Indeed, a cycle $C_k$ or a path $P_k$ is avoidable as one takes a Hamiltonian orientation of $C_k$/$P_k$ and let $f(e)$ be the next edge of that orientation. For a tree $T$ with $\Delta(T)=3$ and exactly one vertex $u$ of degree 3, if the neighbors of $u$ are $v_1,v_2,v_3$, then we can let $f$ be defined as $uv_1\to uv_2 \to uv_3 \to uv_1$ and the edges of the ``hanging paths'' are sent to the edge "one closer" to $u$. So all we are left to show is that all graphs $G\in \cM(P_3)$ with $\Delta(G)=3$ are $C_k^+$ and $B_k$. As some $C_k^+$ is a subgraph of any connected cyclic graph with maximum degree 3, and some $B_k$ is a subgraph of any tree
containing at least two vertices of degree 3, it is enough to show that $B_k$ and $C_k^+$ are unavoidable for $P_3$.

To avoid an $f$-free $P_3$ in $C_k^+$, the edges of the cycle must be mapped cyclically by $f$. So we may assume that $v_2$ is the vertex of degree 3, $f(v_1v_2)=v_2v_3$, $f(v_2v_3)=v_3v_4$ and $u$ is the leaf of $C_k^+$. But then if $f(uv_2)=v_2c_3$, then $v_1v_2u$ is $f$-free, while if $f(uv_2)=v_1v_2$, then $uv_2v_3$ is $f$-free. (If $f(uv_2)$ is any other edge, then both of these $P_3$s are $f$-free.)

Let $u$, $v$ be the degree-3 vertices of $B_k$ and let $u=w_0,w_1,w_2,\dots,w_{k-1},v=w_{k}$ be the vertices of the path from $u$ to $v$. To avoid an $f$-free $P_3$, for every $i=1,2,\dots,\ell$ at least one of the two edges incident to $w_i$ must be sent by $f$ to the other such edge. Also, to avoid an $f$-free $P_3$, the edges incident to $u$ must be mapped by $f$ cyclically, and the same is true for $v$. This shows that $G$ is unavoidable for $P_3$ if $\ell=0$, i.e. $u$ and $v$ are adjacent. Note that if $u$ and $v$ are adjacent and have common neighbors, then $G$ contains $C_3^+$ as a proper subgraph, so $G\notin \cM(P_3)$. If $u$ and $v$ are not adjacent, then because of the cyclical mapping of the edges adjacent to $u$ and to $v$, to avoid $f$-free $P_3$, we must have $f(w_1,w_2)=vw_1=w_0w_1$. And this implies $f(w_3w_2)=w_2w_1$ and so on. But the same reasoning at the other end of the path shows $f(w_iw_{i+1})=w_{i+1}w_{i+2}$, so for some $i$ we must have $f(w_{i-1}w_i)=w_{i-2}w_{i-1}$, $f(w_iw_{i+1})=w_{i+1}w_{i+2}$ and so $w_{i-1},w_i,w_{i+1}$ form an $f$-free copy of $P_3$.
\end{proof}

For the characterization of $\cM(2K_2)$ the following
 observation will be useful.

\begin{lemma}   \label{l:colin}
A graph $G$ admits a $2K_2$-free edge mapping $f\in F_{G,1}$
 if and only if each component in the complement $\overline{L(G)}$
  of its line graph is either a tree or a unicyclic graph.
\end{lemma}

\begin{proof} As noted after Lemma \ref{outdeg}, every mapping $f:E(G)\to E(G)$ with the property $f(e)\ne e$
 for all $e\in E(G)$ can be represented with a digraph $\overrightarrow{\Gamma}_f$
 of constant out-degree 1 on the vertex set $V(\overrightarrow{\Gamma}_f)=E(G)$, by
 drawing an arc from $e$ to $e'$ if $f(e)=e'$.
A pair $e,e'\in E(G)$ of disjoint edges is an
 $f$-free copy of $2K_2$ if and only if neither $f(e)=e'$ nor
  $f(e')=e$ holds.
Hence if no $f$-free $2K_2$ occurs, then all edges of
$\overline{L(G)}$ correspond to an arc in $\overrightarrow{\Gamma}_f$ (where
 cycles of length 2 are allowed).
Here out-degree 1 implies that every component is
 a tree or a unicyclic graph.

Conversely, if each component of $H=\overline{L(G)}$ is a tree
 or a unicyclic graph, then $H$ admits an orientation, say $\overrightarrow{\Gamma}$
 (with exactly one 2-cycle in each tree component by doubling
 an edge of the component in question) to which
 an edge mapping $f$ is naturally associated.
Then $G$ contains no $f$-free copy of $2K_2$ under this $f$,
 because the copies of $2K_2$ in $G$ are in one-to-one
 correspondence with the non-edges of the line graph $L(G)$.
\end{proof}

\begin{proposition}\label{m2k2} ~~~
    $$
      \cM(2K_2)=\{4K_2, P_3+2K_2, P_3+K_3,P_3+K_{1,3},
       C_4+K_2, P_6,W_4^-,C_4^{++}\}
    $$
    where $W_4^-$ is the graph on vertices $a,b,c,d,e$ with $W_4^-[a,b,c,d]=C_4$ and $d(e)=3$, while $C_4^{++}$ has vertex set $a,a',b,c,c',d$ with $C_4^{++}[a,b,c,d]=C_4$, $d(a')=d(c')=1$ and $(aa'),(cc')\in E(C_4^{++})$, $(ac)\notin C_4^{++}$.
\end{proposition}

\begin{proof} We first show that the eight listed graphs are
 members of $\cM(2K_2)$.
We apply Lemma \ref{l:colin} to each graph.
 \begin{itemize} \item $\overline{L(4K_2)}=K_4$ has more than one cycle, and $\overline{L(4K_2-e)}=K_3$ is unicyclic.

\smallskip

\item $\overline{L(P_3+2K_2)})=K_4-e$ has more than one cycle, and $\overline{L(P_3+2K_2-e)}$ is either $K_3$ or $P_3$.

\smallskip

\item $\overline{L(P_3+K_3)}=K_{2,3}$ has more than one cycle, and $\overline{L(P_3+K_3-e)}$ is either $C_4$ or $K_{1,3}$.

\smallskip

\item $\overline{L(P_3+K_{1,3})}$ is the same as the preceding case $\overline{L(P_3+K_3)}$, and also we have the same options when an edge is removed.

\smallskip

\item $\overline{L(C_4+K_2)}=K_1\lor 2K_2$ has more than one cycle, and $\overline{L(C_4+K_2-e)}$ is either $2K_2$ or the paw graph (a triangle with a pendant edge).

\smallskip

\item $\overline{L(P_6)}=\overline{P_5}$ has more than one cycle, and $\overline{L(P_6-e)}$ is either $P_4$ or $C_4$ or the paw.

\smallskip

\item $\overline{L(W_4^-)}$ is the 6-cycle plus one vertex adjacent to two antipodal vertices of $C_6$, also called the theta-graph $\theta_{3,3,2}$. It has more than one cycle, and $\overline{L(W_4^--e)}$ is either $C_6$ or $C_5$ with a pendant edge or $P_5$ with a pendant vertex adjacent to the middle vertex.

\smallskip

\item $\overline{L(C_4^{++})}$ is $C_6$ with a long diagonal. It has more than one cycle, and $\overline{L(C_4^{++}-e)}$ is either $P_5$ or $C_4$ with a pendant edge. \end{itemize}

In proving that $\cM(2K_2)$ contains no other graphs, we shall
 use the fact that all those eight graphs above do belong to $\cM(2K_2)$.
Several case distinctions will be needed.
In the sequel $G$ denotes any supposed member of $\cM(2K_2)$,
 and the cases will end up either with the contradiction
 $G\notin \cM(2K_2)$ or with the conclusion that $G$ is
 an already known member of $\cM(2K_2)$.

\newcommand{\cass}[1]{\medskip \noindent\textsl{Case #1:\/} }

\medskip \noindent
\underline{First, assume that $G$ is disconnected.}

\cass{1.1} $G$ has at least three components.

In case of four components or more, we can take
 one edge from each and obtain
 $4K_2\in \cM(2K_2)$, hence $G=4K_2$.
If $G$ has three components, then
 at least one of them contains $P_3$ because $G\neq 3K_2$.
Taking a $P_3$ from the largest component and one edge from each
 of the other two components we obtain $P_3+2K_2\in \cM(2K_2)$,
  hence $G$ cannot have further edges.

\cass{1.2} $G$ has two components.

Assume that $G=G_1+G_2$, where $|E(G_1)|\leq |E(G_2)|$.
If $G_i$ has $m_i$ edges for $i=1,2$ then $G$ contains at least
 $m_1m_2$ copies of $2K_2$, while the number of edges is $m_1+m_2$.
We cannot have $m_1>2$ (what would imply $m_2\geq 3$)
 for otherwise keeping only two edges of
 $G_1$ the remaining graph would still contain at least $2m_2$
 copies of $2K_2$ while there would be only $m_2+2$ edges, i.e.\
 fewer edges than $2K_2$ subgraphs, contradicting the assumption
 that $G$ is a (minimal) member of $\cM(2K_2)$.
Consequently $G_1=P_3$ or $G=K_2$ holds.

In the same way, we see that if $G=P_3$, then $m_2\leq 3$,
 because for a larger $m_2$, we may omit an edge from $G_2$.
Note that $2P_3\notin \cM(2K_2)$, we have already seen this fact
 during the verification of $P_3+K_3\in \cM(2K_2)$.
Hence if $G_1=P_3$, then $m_2=3$ and $G_2$ is either $K_3$ or
 $K_{1,3}$ or $P_4$.
The last one is excluded, however, since $P_3+2K_2$ is a proper
 subgraph of $P_3+P_4$.

It remains to consider $G_1=K_2$.
In this case we have $\overline{L(G)}=K_1\lor \overline{L(G_2)}$,
 always connected.
This $\overline{L(G)}$ is a star if $G_2$ contains no $2K_2$, unicyclic if
 $G_2$ has precisely one $2K_2$ subgraph, and has more than
 one cycle if $G_2$ contains more than one copy of $2K_2$.
Hence the latter should hold for $G_2$.
Also, $P_3+K_2\not\subset G_2$, for otherwise we would have
 $P_3+2K_2\subset G$.
Moreover, $G_2$ contains a non-star spanning tree.
It follows that $|V(G_2)|=4$, and two copies of $2K_2$ yield
 $C_4\subset G_2$; in fact $G_2=C_4$, due to the minimality
 condition on the members of $\cM(2K_2)$.

\medskip \noindent
\underline{Second, assume that $G$ is a tree.}

\def \diam {\mathrm{diam}(G)}

\cass{2.1} $\diam\leq 3$.

In the case of $\diam\leq 2$, $G$ would be
 a star, not containing any $2K_2$ subgraph.
On the other hand $\diam=3$ would mean that $G$ is a
 double star and its middle edge would not be contained
 in any $2K_2$, contradicting the minimality condition
 concerning the members of $\cM(2K_2)$.

\cass{2.2} $\diam\ge 4$.

We can immediately settle the case $\diam\geq 5$, because then
 $G$ contains $P_6$, hence it must be $P_6$ if $G\in \cM(2K_2)$.

Assume $\diam=4$ and let $P=v_1v_2v_3v_4v_5$ be a longest path in $G$.
Any further vertices must be leaves adjacent to $P$, otherwise
 a path of length 2 attached to $P$ would yield $P_3+2K_2$ as a
  proper subgraph (and also a path longer than $P$, if attached
   to $v_2$ or $v_4$).
The same would occur if $v_3$ had two further neighbors, or if
 both $v_2$ and $v_3$ (or both $v_3$ and $v_4$) had a further neighbor.
But if just one new leaf is attached to $v_3$, then $\overline{L(G)}$
 is the unicyclic ``bull'' graph, and if it is attached to $v_2$
 or $v_4$ then $\overline{L(G)}$ is the 4-cycle with a pendant edge.
On the other hand, if more than one new leaf is attached to
 $v_2$ or $v_4$, then a proper subgraph $P_3+K_{1,3}\in \cM(2K_2)$ occurs.
The same happens if both $v_2$ and $v_4$ get a new leaf neighbor.

\medskip \noindent
\underline{Third, assume that $G$ is connected and contains a cycle.}

\cass{3.1} The longest cycle has a length of at least 5.

Cycle lengths of at least 6 are excluded by the fact $P_6\in \cM(2K_2)$.
Hence assume that $C$ is a 5-cycle in $G$.
Any neighbor of $C$ would yield a $P_6$ subgraph, consequently
 $|V(G)|=5$.
If there are two disjoint chords in $C$, then we obtain
 $W_4^-\subset G$, and hence $G=W_4^-$ follows.
On the other hand, if $C$ does not have disjoint chords, then
 $\overline{L(G)}$ is a 5-cycle supplemented with at most two
 pendant edges, hence unicyclic and therefore, by Lemma \ref{l:colin}, $G$ cannot be a
 member of $\cM(2K_2)$. 

\cass{3.2} The longest cycle has length 4.

Let $C=v_1v_2v_3v_4\subset G$.
Then the set $Z:=V(G)\setminus V(C)$ is independent, otherwise
 $G=C_4+K_2$ would follow.
Suppose first that there are two disjoint edges from $C$ to $Z$.
Their ends cannot be consecutive along $C$, because then $P_6$
 would be a proper subgraph of $G$.
Hence their ends in $C$ are antipodal, thus $C_4^{++}\subset G$
 holds and $G=C_4^{++}$ follows.

The other possibility is that all edges between $C$ and $Z$
 share a vertex, which can be in $C$ or in $Z$.

If it is $v_i\in V(C)$, say $v_1$, then it has $|Z|$
 pendant neighbors.
For $|Z|\geq 3$ a subgraph $P_3+K_{1,3}\subset G$ would occur,
 that cannot be the case.
For $|Z|=2$ the presence of the edge $v_2v_4$ would also lead to
 the contradiction $P_3+K_3\subset G$, moreover the edge
  $v_1v_3$ would not be included in any copy of $2K_2$,
  hence it cannot be present.
So in that case the only possibility for $G$ remains to be the
 4-cycle with two pendant edges at $v_1$.
But then $\overline{L(G)}$ is the unicyclic graph formed by a
 $C_4$ with two pendant edges at antipodal vertices, hence,  by Lemma \ref{l:colin},
 such a $G$ cannot belong to $\cM(2K_2)$.

Finally if $|Z|=1$, then either there is just one connecting
 edge say $v_1z$, or two edges $v_1z$ and $v_3z$, because $G$
 does not contain a 5-cycle.
In the first case $V(C)$ may induce $K_4$, but even in that
 densest situation $\overline{L(G)}$ is a tree (three
 paths of length 2 at a vertex of degree 3 as their common end).
In the second case $G$ would be $K_{2,3}$ and hence
 $\overline{L(G)}=C_6$ would hold, consequently $G\notin \cM(2K_2)$.

\cass{3.3} Every cycle in $G$ is a triangle.

Let $T=v_1v_2v_3\subset G$ be a triangle.
Since $P_3+K_3$ and $P_3+2K_2$ are in $\cM(2K_2)$, the set
 $V(G)\setminus V(T)$ either is independent or induces just
 one edge.

In the first case, every vertex outside $T$ is a pendant vertex
 with exactly one neighbor in $T$, because $G$ is connected and
 does not contain a 4-cycle.
Let $X_i$ denote the set of neighbors of $v_i$ outside $T$,
 for $i=1,2,3$.
None of the $X_i$ is empty, because e.g.\ if $X_3=\varnothing$,
 then the edge $v_1v_2$ does not occur in any $2K_2$.
On the other hand, if some $X_i$ has more than one vertex, then
 we get the contradiction $P_3+2K_2\subset G$.
Hence $|X_1|=|X_2|=|X_3|=1$ should hold.
But then $\overline{L(G)}\cong G$ would be a unicyclic graph.

In the other case, let $e=v_4v_5$ be the unique edge outside $T$,
 and let $v_3v_4$ be an edge connecting $e$ with $T$.
Then $v_3v_5$ may be an edge or a non-edge.

If $v_3v_5$ is an edge, then a pendant edge at a $v_i\neq v_3$
 would yield $P_3+K_3\subset G$; and in case of more than one
 pendant edge at $v_3$ we would have $P_3+2K_2\subset G$.
On the other hand, if there is only one pendant edge at $v_3$,
 then $\overline{L(G)}$ is a unicyclic graph, namely a triangle
 with two vertices having two pendant neighbors each, hence
 $G\notin \cM(2K_2)$.

In the last case, if $v_3v_5$ is not an edge, there can be a
 pendant edge at $v_3$ or at $v_1$ (or $v_2$) but not both,
 otherwise, we would have $P_3+2K_2 \subset G$.
For the same reason a second pendant neighbor at the same vertex
 is not possible.
It follows that $\overline{L(G)}$ is the unicyclic graph formed
 by a triangle with two pendant neighbors at one vertex and
 two pendant neighbors at another vertex.

\medskip

The above analysis shows that $\cM(2K_2)$ contains no graphs
 beyond the eight listed ones.
\end{proof}

Now we start working toward Theorem \ref{forest}.

\begin{lemma}\label{3t-2}
For any $t\ge 2$, we have $(3t-2)K_2\in \cM(tK_2)$.    
\end{lemma}

\begin{proof}
    First we show that every $f \in F_{(3t-2)K_2,1}$ admits an $f$-free $tK_2$.
Let $e_1,...,e_m$ be the edges of $(3t-2)K_2$.
For any given $f$, construct the auxiliary digraph $F$ with vertex set $u_1,...,u_m$
where there is an arc from $u_i$ to $u_j$ if and only if $f(e_i) = e_j$.
Every vertex has out-degree 1, therefore each component in the undirected
underlying graph $F'$ of $F$ is a tree or a unicyclic graph.
Hence $F'$ is 3-colorable and therefore it contains an independent set of size $t$.
This set corresponds to an $f$-free $tK_2$ in $(3t-2)K_2$.
 
To show that $G = mK_2$ with $m=3t-3$ admits an $f$ without $f$-free $tK_2$ subgraphs,
take the mapping $e_{3i-2} \to e_{3i-1} \to e_{3i}\to  e_{3i-2}$ for $i = 1,2,...,t-1$.
\end{proof}

\begin{proposition}\label{matching}
    The family  $\cM(tK_2)$  is finite for every $t \geq 2$.
\end{proposition}

\begin{proof}
    Let $G = (V,E)$ be any graph from $\cM(tK_2)$, and
let $M=\{e_1,...,e_m\}$ be the edge set of a largest matching in $G$.
We certainly have $m \geq t$, otherwise $G$ would not be unavoidable; and by the
minimality of graphs in $\cM(tK_2)$, we also know from  Lemma \ref{3t-2} that $m \leq 3t-2$.
(In fact if $m=3t-2$, then $G = mK_2$, but this is unimportant.)
We need to prove that $|V|$ is bounded above by a function of $t$.
We have seen that $V' := V(M)$ is bounded by $|V'| = 2m < 6t$.
 
Since the matching $M$ is non-extendable, $V'':=V\setminus V'$ is an independent set.
We label each $x \in V''$ with its neighborhood $N(x)$, which is a subset of $V'$.
This classifies the vertices of $V''$ into $2^{2m} - 1$ types $T_1,T_2,...,T_k$.
The distribution of types is represented by the $k$-tuple $(n_1,...,n_k)$, whose
meaning is that exactly $n_i$ vertices of $V''$ have type $T_i$.
In this way each $G \in \cM(tK_2)$ is represented by
\begin{enumerate}
    \item 
    the subgraph induced by $V'$,
    \item
    the associated $k$-tuple $(n_1,...,n_k)$.
\end{enumerate}
Since all graphs in $\cM(tK_2)$ are minimal under inclusion, it follows that
if two graphs have the same subgraph under (1), then their associated $k$-tuples
are incomparable, i.e., each of them has a component $n_i$ larger than the same
component in the representation of the other graph.
 
Assume now for a contradiction that $\cM(tK_2)$ is infinite.
We are going to select a sequence from its members as follows.
Since $m < 6t$, an infinite subfamily $\cM_1(tK_2)$ has the same subgraph under (1).
Restricting attention to $\cM_1(tK_2)$, select an arbitrary member $G_1\in \cM_1(tK_2)$.
Let $i_1$ be its component where $n_{i_1}$ exceeds the corresponding component
of infinitely many members of $\cM_1(tK_2)$.
All those members have this component smaller than $n_{i_1}$, hence infinitely many
of them have the same value in the $i_1$-component.
Let $\cM_2(tK_2) \subset \cM_1(tK_2)$ be the family of the latter graphs.
Analogously, we select an arbitrary $G_2$ from $\cM_2(tK_2)$.
It has a component indexed with $i_2$ where $n_{i_2}$ exceeds the corresponding
component of infinitely many members of $\cM_2(tK_2)$, among which infinitely many
have the same value in the $i_2$-component.
Note that $i_1 = i_2$ is impossible, by the definition of $\cM_2(tK_2)$.
After $k$ steps of this kind we obtain an infinite subfamily $\cM_k(tK_2)$ of $\cM(tK_2)$
in which all graphs have the same associated $k$-tuple $(n_1,...,n_k)$.
However, with a fixed $k$-tuple and with a fixed subgraph induced by $V'$ we specify
just one graph.
This contradiction proves the proposition.
\end{proof}

The other direction of the proof of Theorem \ref{forest} will be shown in two cases. The next proposition considers star-forests.

\begin{proposition}\label{treeforest}
    \
$\cM(G)$ is infinite for all star-forests that are not matchings, i.e.\ whenever $G=\bigcup_{i=1}^ja_iK_{1,r_i}$ for some $r_1\ge r_2 \ge \dots \ge r_j\ge 1$ with $r_1\ge 2$ and the $a_i$ are positive integers.
\end{proposition}

\begin{proof}
We define the graph class $D_{\ell,k,\ell'}$ as follows: we start with  the double-star $D_{\ell,\ell}$ with center $\{u,v\}$. Then we subdivide the edge $uv$ with $k-1$ intermediate vertices $w_1,w_2,\dots,w_{k-1}$ and add $\ell'$ pendant edges to all $w_i$. First, we claim that $D=D_{3r-4,k,3r-6}\in \cM(K_{1,r})$ holds for all $r\ge 2$, $k\ge 1$. Let us fix $f\in F_{D,1}$. We need to show that there exists an $f$-free copy of $K_{1,r}$. Observe that $u$ and $v$ have degree $3r-3$ in $D$. Applying Lemma \ref{outdeg} (2), we obtain that the only way there is no $f$-free copy of $K_{1,r}$ with center $u$ or $v$ is that edges incident to $u$ are partitioned into triples $(e,e',e'')$ such that $f(e)=e',f(e')=e'',f(e'')=e$. In particular, the edge $uw_1$ is mapped to an edge not incident to $w_1$. Therefore if $f(w_1w_2)$ is not incident to $w_1$, then applying Lemma \ref{outdeg} (2) with $m=2$ and $n-m=3r-6$ to the edges incident to $w_1$, we would have an $f$-free copy of $K_{1,2+\lceil \frac{3r-8}{3}\rceil}=K_{1,r}$. So $f(w_1w_2)$ must be incident to $w_1$ and thus non-incident to $w_2$. Repeating this argument shows that every $f(w_iw_{i+1})$ should be incident to $w_i$ and not to $w_{i+1}$ or there is an $f$-free copy of $K_{1,r}$. But we could have started this reasoning at $v$ and obtain the same way that $f(w_{k-1}v)$ is incident to $v$ and thus not to $w_{k-1}$ and then all $f(w_iw_{i+1})$ should be incident to $w_{i+1}$ and not to $w_i$. This contradiction shows that there must exist an $f$-free copy of $K_{1,r}$, i.e.\ $D$ is unavoidable for $K_{1,r}$.

We also need to show that $D$ is \textit{minimal} unavoidable for $K_{1,r}$; that is, $D\setminus e$ is not unavoidable for $K_{1,r}$ for any $e\in E(D)$. We distinguish three types of edges of $D$: $e$ is a \textit{path-edge} if $e$ lies on the path from $u$ to $v$, $e$ is an \textit{end-leaf edge} if $e$ is incident to $u$ or $v$ and not a path-edge, and $e$ is a \textit{path-leaf edge} if $e$ is incident to some $w_i$ and not a path-edge. In all three cases, we need to define an $f\in F_{D\setminus e,1}$ that admits no $f$-free copy of $K_r$.

\begin{itemize}

    \medskip 
    
    \item 
    If $e$ is an end-leaf edge, then by symmetry we may assume $e$ is incident to $u$. There are $3r-5$ other end-leaf edges incident to $u$ which we partition into $r-2$ triples which $f$ maps cyclically and the last end-leaf edge incident to $u$ is mapped to $uw_1$. Therefore the largest $f$-free star centered at $u$ is of size $r-1$. The $3r-6$ path-leaf edges incident to $w_i$ are partitioned in $3r-2$ triples that are mapped cyclically by $f$ and we let $f(uw_1)=w_1w_2$, $f(w_iw_{i+1})=w_{i+1}w_{i+2}$ for all $i=1,2,\dots,j-1$. This ensures that the largest $f$-free star centered at any $w_i$ is of size $r-1$. Finally, $w_jv$ and the end-leaf edges incident to $v$ are partitioned into triples that are mapped cyclically, so the largest $f$-free star centered at $v$ is also of size $r-1$.
    
    \medskip 
    
    \item
    If $e$ is a path-leaf edge adjacent to $w_i$, then let $f(w_hw_{h+1})=w_{h-1}w_h$ if $h+1\le i$, and let $f(w_hw_{h+1})=w_{h+1}w_{h+2}$ if $h\ge i$. The $3r-3$ edges incident to $u$ are mapped cyclically by $f$ in $r-1$ triples, similarly to those that incident to $v$, and those path-leaf edges that are incident to $w_h$ for some $h\neq i$. Finally, let $3r-9$ out of the remaining path-leaf edges adjacent to $w_i$ be partitioned into $r-3$ triples and let $f$ map them cyclically, and for the two remaining such edges $e',e''$ let $f(e')=w_{i-1}w_i$ and 
$f(e'')=w_iw_{i+1}$. It is easy to check that $f$ does not admit any $f$-free copies of $K_{1,r}$.

\item
Finally, if $e$ is a path-edge, then the removal of $e$ yields two tree components, each being a proper subgraph of some graph already discussed under end-leaf deletion.
\end{itemize}

Now we consider the general case $G=\bigcup_{i=1}^ja_iK_{1,r_i}$ with $r_1\ge r_2\ge \dots \ge r_j\ge 1$ and $r_1\ge 2$. Let $a=\sum_{i=1}^ja_i$ and let $G_{a,k}$ consist of $(2r_1+1)a$ copies of $D_{3r_1-4,k,3r_1-6}$. 

\begin{clm}
    $G_{a,k}$ is unavoidable for $G$.
\end{clm}

\begin{proof}[Proof of Claim]
For $f\in F_{G,1}$ and $1\le h \le (2r_1+1)a$, let $f_h$ defined on the $h$th copy of $D_{3r_1-4,k,3r_1-6}$ by $f_h(e)=f(e)$, whenever $f(e)$ is in the same copy of $D_{3r_1-4,k,3r_1-6}$, and otherwise $f_h(e)$ is picked arbitrarily with $f_h(e)\neq e$. By the first part of the proof, $D_{r_1-4,k,r_1-6}$ is unavoidable for $K_{1,r_1}$, and so for all $h$ there exists an $f_h$-free copy $K^h$ of $K_{1,r_1}$ in the $h$th copy of $D_{r_1-4,k,r_1-6}$. These are clearly, $f$-free as well, but it is not guaranteed that their union (or the union of some of them) is still $f$-free as for some $e\in K^h,e'\in K^{h'}$ we might have $f(e)=e'$. 

With the notation after Lemma \ref{outdeg} us, let introduce the directed graph $\overrightarrow{\Gamma}_{f,\cE}$ such that $\cE$ is the partition with parts being the edge sets of distinct $K^h$s.
 As $K^h$ has $r_1$ edges, the maximum out-degree of $\overrightarrow{\Gamma}_{f,\cE}$ is at most $r_1$, so by Lemma \ref{outdeg} (1), the underlying undirected graph contains an independent set $I$ of size $\frac{(2r_1+1)a}{2r_1+1}=a$. Clearly, $\bigcup_{h\in I}K^h$ is an $f$-free copy of $aK_{1,r}$ that contains $G$.
\end{proof}

$G_{a,k}$ is unlikely to be \textit{minimal} unavoidable for $G$, but if we remove edges of $G_{a,k}$ greedily according to whether the remaining graph is unavoidable or not, then we obtain $G'_{a,k}\in\cM(G)$. Observe that independently of the greedy order of the edges, $G'_{a,k}$ must contain at least one (actually many more) complete copy of $D_{3r_1-4,k,3r_1-6}$. If not, then by the fact that $D_{3r_1-4,k,3r_1-6}\in \cM(K_{1,r})$, there would exist a mapping $f_h$ for the remaining part of the $h$th copy of $D_{3r_1-4,k,3r_1-6}$ that does not admit a free $K_{1,r_1}$. Considering the union $f$ of these $f_h$s, $f$ would not admit a free copy of $K_{1,r_1}$ let alone that of $G$.

Finally, observe that for $k\neq k'$ none of the graphs $D_{3r_1-4,k,3r_1-6}$ and $D_{3r_1-4,k',3r_1-6}$ contains the other, so neither does $G'_{a,k}$ and $G'_{a,k'}$. So $\{G'_{a,k}:k\in \mathbb{N}\}$ is an infinite subset of $\cM(G)$.
\end{proof}

The next proposition is a general statement that resolves the case of all non-star forests in Theorem \ref{forest}.

\begin{proposition}\label{p4avoid}
    Let  $\cH$ be a family of graphs with $\ex(n,\cH)\le c_\cH n$ for some constant $c_\cH$ such that every $H\in \cH$ contains a $P_4$. Then $\cM(\cH)$ is infinite.
\end{proposition}

\begin{proof}
    We need several auxiliary results.

    \begin{clm}\label{treep4}
        Every tree is avoidable for $P_4$ and thus for any graph that contains a $P_4$.
    \end{clm}

    \begin{proof}
        Let $T$ be an arbitrary tree and $r$  an arbitrarily chosen root in $T$. Let us define $f\in F_{T,1}$ by letting $f(e)=e'$ for any $e$ not incident with $r$ with $e'$ being the edge one closer to $r$ and adjacent to $e$. We define $f$ on edges $e$ incident with $r$ arbitrarily but with $f(e)\neq e$. Observe that for any $P_4$ in $T$ if $v$ is the vertex closest to $r$, then an edge $e$ of $P_4$ that is not incident to $v$ is $f$-mapped to another edge of this $P_4$. Therefore there is no $f$-free copy of $P_4$ in $T$.
    \end{proof}

    \begin{theorem}[Erd\H os, Sachs \cite{ESa}]\label{esa}
    For any $d$ and $k$, there exists a $d$-regular graph with girth $k$.
    \end{theorem}

    \begin{clm}\label{nonbip}
    Any $d$-regular graph $G$ with $d\ge 6c_\cH$ is unavoidable for $\cH$.
    \end{clm}

    \begin{proof}
    Consider any $f\in F_{G,1}$. Applying Lemma \ref{outdeg} (1), we obtain a subgraph $G'$ of $G$ such that the average degree of $G'$ is at least $d/3$ and for any $e\in E(G')$ we have $f(e)\in E(G)\setminus E(G')$. Clearly, any subgraph of $G'$ is $f$-free, and since the average degree of $G'$ is at least $2c_\cH$, $G'$ contains some $H\in \cH$.
    \end{proof}

    With all the above statements in hand, we are ready to prove the theorem. Applying Theorem \ref{esa}, we take an infinite sequence of $6c_\cH$-regular graphs, one graph $G_k$ for each possible girth $k$. This $G_k$ is unavoidable for $\cH$, by Claim \ref{nonbip}. So there exists $G'_k \subseteq G_k$ with $G'_k\in \cM(\cH)$. This $G'_k$ is not a tree, by Claim \ref{treep4}, i.e. $G'_k$ contains a cycle $C_k$. We can begin with $k_1=5$. Then choose $k_2 > |V(G'_{k_1})|$. Clearly, $G'_{k_1}$ is not a subgraph of $G'_{k_2}$ because the latter has higher girth, and $G'_{k_1}$ does not contain $G'_{k_2}$ as a subgraph because the latter has more vertices. So an infinite sequence of mutually incomparable members of $\cM(\cH)$ can be selected.
\end{proof}

Choosing $\cH$ to be $\{F\}$ for some non-star forest, $\cC_k$ the cycles of length at least $k$, $\cM K_p$ the family of $K_p$-minor-free graphs, or $\cS K_p$ the family of topological $K_p$-minor-free graphs, one obtains that the corresponding set of minimal unavoidable graphs is infinite.

\medskip

Now we have all the ingredients for the proof of Theorem \ref{forest}.

\begin{proof}[\textbf{Proof of Theorem \ref{forest}}]
 If $F$ is a matching, then $\cM(F)$ is finite by Proposition \ref{matching}.
If $F$ is a star-forest that is not a matching, then $\cM(F)$ is infinite by Proposition \ref{treeforest}.
If $F$ is a non-star forest, then $F$ contains a $P_4$, and therefore $\cM(F)$ is infinite by Proposition \ref{p4avoid} applied with $\cH=\{F\}$.
\end{proof}

In the remainder of this section, we show that $\cM(G)$ is infinite for some other graphs $G$.

\begin{proposition}\label{mk3}
    $\cM(K_3)$ is infinite.
\end{proposition}

\begin{proof}
    Observe that for any graph $G$, edge $e\in E(G)$ and $f\in F_{G,1}$ there is at most one triangle $T$ such that $T$ is not $f$-free because $T$ contains both $e,f(e)$. Therefore,  any graph $G$ is unavoidable that has more triangles than edges. In particular, $K_6$ minus an edge is unavoidable. Let us introduce two of its subgraphs: $H_1:=K_6\setminus E(P_3)$, with 13 edges and 13 triangles, and $H_2:=K_6\setminus E(2K_2)$, with 13 edges and 12 triangles. We call the edges of $K_6\setminus H_i$ the \textit{missing edges of $H_i$}.
    We construct $G_k$ as a ``chain'' of $k-2$ copies of $H_2$, any two consecutive copies sharing an edge, and a copy of $H_1$ is taken at each end of the chain again sharing an edge with the first and last copy of $H_2$ in the chain. So, $G_k$ has $4k+2$ vertices and $12k+1$ edges. The overall number of triangles is $12k+2$, and thus $G_k$ is unavoidable for all $k$. Note that $G_k$ is not uniquely determined at this moment, since we have not specified which edges of neighboring parts are shared. To prove that with appropriately chosen shared edges $G_k$ belongs to $\cM(K_3)$, we will need the following claims.

    \begin{clm}\label{h2}
        Let $e_1$ be the edge of $H_2$ that extends the missing edges of $H_2$ to a perfect matching, and let $e_2$ be an edge of $H_2$ that connects the two missing edges of $H_2$. Then for any $e\in E(H_2)$ (including the possibility of $e=e_1$ or $e=e_2$), there exists an injection $\iota$ from the triangles of $H_2\setminus e$ to $E(H_2)\setminus \{e,e_1,e_2\}$ with $\iota(T)\subseteq T$.
    \end{clm}

\begin{proof}
    First observe that for any $e\in E(H_2)$, the number of triangles in $H_2\setminus e$ is at most 10. We will use the shadow function $KK(m)$ that tells us the minimum number of pairs that are contained in any set of $m$ triples. By the famous Kruskal-Katona theorem, its value is attained at the initial segment of the colex order of triples. So for $m=1,2,3,4,5,6,7,8,9,10$, the value of $KK(m)$ is $3,5,6,6,8,9,9,10,10,10$, respectively.

    We need to verify that Hall's condition holds in the auxiliary bipartite graph with one part $E(H_2)\setminus \{e,e_1,e_2\}$ and other parts of the triangles of $E(H_2)\setminus e$, and $e^*$ joined to $T$ if and only if $e^*\subset T$. As only $e_1,e_2$ can make a problem, a set of vertices corresponding to $m$ triangles are joined to at least $KK(m)-2$ vertices corresponding to edges. So if $m\le 8$, then $KK(m)-2\ge m$ ensures Hall's condition. The case $m=10$ is trivial as all edges are contained in at least one triangle. Finally, the case $m=9$ follows from the fact that every edge is contained in at least 2 triangles, so the neighborhood of vertices corresponding to 9 triangles is already the whole ``edge part'' of the auxiliary graph.
\end{proof}

    \begin{clm}\label{h1}
Let $e'_1$ be an edge of $H_1$ that is not adjacent to any of the missing edges of $H_1$. Then for any $e\in E(H_1)$ (including the possibility of $e=e'_1$), there exists an injection $\iota$ from the triangles of $H_1\setminus e$ to $E(H_1)\setminus \{e,e'_1\}$ with $\iota(T)\subseteq T$.
    \end{clm}

    \begin{proof}
        The proof is analogous to that of Claim \ref{h2}. Here $H_1\setminus e$ has at most 11 triangles. But since this time only $e'_1$ can cause problems, $m$ triangles have a neighborhood of size at least $KK(m)-1$, so this covers the case $m\le 9$. The case $m=11$ is trivial, and case $m=10$ is analogous to the case $m=9$ of Claim \ref{h2}. 
    \end{proof}
Observe that the claims above show the existence of mappings $f$ with no $f$-free $K_3$ in $H_j\setminus e$: if $e'=\iota(T)$ for some triangle $T$, then define $f(e')$ to be one of the other edges of $T$. Because of this, we define $G_k$ in such a way that a copy of $H_2$ shares $e_2$ with its right neighbor piece in the chain, and $e_1$ with its left neighbor piece in the chain. The copies of $H_1$ share $e'_1$ with their neighbor piece in the chain.

Let $e^*$ be an arbitrary edge of $G_k$ and consider $G_k\setminus e^*$. If $e^*$ belongs to the $h$th (and $(h+1)$st) piece(s) of the chain, then we use Claims \ref{h2} and \ref{h1} to define $f$ on this (these) piece(s). After that, an edge shared by two neighbor pieces is used in the piece that is farther from the piece from which $e^*$ has been removed. In this way, an $f$-free mapping is obtained for every proper subgraph of $G_k$.
\end{proof}

\begin{proposition}
    Let $G_0$ be a graph with an infinite $\cM(G_0)$ and let $G_1,\dots,G_k$ be distinct  subgraphs  of $G_0$.  Then for $G = \cup_{j=0}^k a_j G_j$, $j = 0,\dots,k$  with $a_0 \ge 1$ the family $\cM(F)$ is infinite.
\end{proposition}

\begin{proof}
    Set $z:= \sum_{j=0}^ka_j$.
Similarly to the case of star-forests, we claim that for any member $H$ of $\cM(G_0)$ the graph $z(2e(G_0) +1)H$  is unavoidable for $F$. To see this consider an $f\in F_{G,1}$. Each copy of $H$ must contain an $f$-free copy of $G_0$. Then we apply Lemma  \ref{outdeg} to the graph $\overrightarrow{\Gamma}_{f,\cE}$, where $\cE={E_1,E_2,\dots,E_{z'}}$ with $z':=z(2e(G_0)+1)$ and $E_i$ is the edge set of the $f$-free copy of $G_0$ in the $i$th copy of $H$. We obtain that the chromatic number is at most $2e(G_0)  +1$, hence there exists an independent set of size at least $z$ in $\Gamma_{f,\cE}$.  The corresponding copy of $zG_0$ is $f$-free, hence it contains a $G$ that is $f$-free.
    
Now $z(2e(G_0) +1)H$ contains a member $H'$ of $\cM(G)$ for every $H\in \cM(G_0)$. This $H'$ must contain at least one copy of $H$ as otherwise, by the minimality of $H$ for $G_0$, we could obtain an $f$ without a free copy of $G_0$ and thus without a free copy of $G_0$. As  any $H_1,H_2\in \cM(G_0)$ are incomparable, so are $H'_1,H'_2$.
\end{proof}

By Proposition \ref{mk3}, $\cM(K_3)$ is infinite, and so we immediately obtain the following corollary.

\begin{corollary}
    $\cM(aK_3 \cup  bP_3 \cup cK_2)$  is infinite for $a>0$.
\end{corollary}

\section{Concluding remarks and open problems}\label{comm}

Here we list those open problems, apart from Conjecture \ref{mconj}, that we find the most natural to consider, and we add some comments regarding the parameters that we addressed in the paper. 

First, let us mention that a relatively standard argument using greedy deletion of $f$-free or $f$-exclusive copies of $G$ and then applying Lemma 
\ref{outdeg} yields the following statement of which the Turán number analog is well-known. We leave the details to the Reader.

\begin{proposition}
    For any graph $G$ and integer $t\ge 2$ there exist constants $c_1(t,G)$ and $c_2(t,G)$ such that
    \begin{enumerate}
        \item 
         $h(n,tG ) \le h(n,G) + c_1(t,G)n$,
         \item
         $s(n,tG)  \le s(n,G) + c_2(t,G)n$. 
    \end{enumerate}
\end{proposition}

\noindent \textbf{Unavoidability for $s(n,G)$}. The observation that made possible the introduction of $\cM(G)$ was that if $H$ is unavoidable for $G$, then so is any supergraph of $H$. Then Proposition \ref{unavoid} proved that $h(n,G)=\ex(n,\cM(G))$ holds. We would like to obtain a similar equality for $s(n,G)$ via minimal $f$-exclusive unavoidable graphs $G$, but the situation is not that simple.

Consider  $C_4$ with vertices $a, b, c, d$ and to each  vertex attach leaf 1 to $a$, 2 to $b$, 3 to $c$, and 4 to $d$, to obtain graph $H$.
The edge mapping $f$ that maps $(1,a)$ to $(b,c)$, $(2,b)$ to $(c,d)$, $(3,c)$ to $(d,a)$, and $(4,d)$ to $(a,b)$ and vice versa belongs to $F_{H,0}$ with no exclusive copy of $2K_2$,  However if we consider $H' = 4K_2\subseteq H$ alone, it is unavoidable for exclusive copies of $2K_2$.


 \begin{problem}
    Determine $h(n, T )$ at least for small trees. In particular, is $3n- 9$ the true value of $h(n, P_4)$?
\end{problem}    

\begin{problem}
    Determine or get good bounds for $s(n,G)$ in case $G$ is $K_{1,r}$, $tK_2$, a tree in general, or a bipartite graph.
\end{problem}

\begin{problem}
    We introduced $q(n, G)$ as an analog of $h(n, G)$, and Proposition \ref{qh} states that these quantities coincide. Similarly, one can define $$p(n, G) := \max\{ m :  \exists f\in F_{n,m,0},
\text{~with no~}
\text{$f$-exclusive copy of~} G ~\text{in}\ K_n  \}.$$ Prove or disprove $p(n,G) = s(n,G)$.    
 \end{problem}

\end{document}